\documentclass[12pt]{amsart}
\usepackage{amsmath}
\usepackage{amssymb}
\usepackage{amsfonts}
\usepackage{latexsym}
\usepackage{amscd}
\usepackage[mathscr]{euscript}
\usepackage{enumitem}
\usepackage{xy} \xyoption{all}
\usepackage[active]{srcltx}

\newcommand{\N}{{\mathbb{N}}}
\newcommand{\Z}{{\mathbb{Z}}}

\newcommand{\uloopr}[1]{\ar@'{@+{[0,0]+(-4,5)}@+{[0,0]+(0,10)}@+{[0,0] +(4,5)}}^{#1}}
\newcommand{\uloopd}[1]{\ar@'{@+{[0,0]+(5,4)}@+{[0,0]+(10,0)}@+{[0,0]+ (5,-4)}}^{#1}}
\newcommand{\dloopr}[1]{\ar@'{@+{[0,0]+(-4,-5)}@+{[0,0]+(0,-10)}@+{[0, 0]+(4,-5)}}_{#1}}
\newcommand{\dloopd}[1]{\ar@'{@+{[0,0]+(-5,4)}@+{[0,0]+(-10,0)}@+{[0,0 ]+(-5,-4)}}_{#1}}
\newcommand{\luloop}[1]{\ar@'{@+{[0,0]+(-8,2)}@+{[0,0]+(-10,10)}@+{[0, 0]+(2,2)}}^{#1}}

\newtheorem{lem}{Lemma}[section]
\newtheorem{corol}[lem]{Corollary}
\newtheorem{theor}[lem]{Theorem}
\newtheorem{prop}[lem]{Proposition}
\theoremstyle{definition}
\newtheorem{defi}[lem]{Definition}

\newtheorem{rema}[lem]{Remark}

\begin{document}

\title[Isomorphisms of Higman-Thompson groups]{The Isomorphism Problem for Higman-Thompson groups}%
\author{Enrique Pardo}
\address{Departamento de Matem\'aticas, Facultad de Ciencias\\ Universidad de C\'adiz, Campus de
Puerto Real\\ 11510 Puerto Real (C\'adiz)\\ Spain.}
\email{enrique.pardo@uca.es}\urladdr{https://sites.google.com/a/gm.uca.es/enrique-pardo-s-home-page/}
\thanks{The author was partially supported by the DGI and European Regional Development Fund, jointly, through Project MTM2008-06201-C02-02, by the Consolider Ingenio ``Mathematica" project CSD2006-32 of the MEC, by PAI III grants FQM-298 and P07-FQM-02467 of the Junta de Andaluc\'{\i}a and by 2009 SGR 1389 grant of the Comissionat per Universitats i Recerca de la Generalitat de Catalunya.} 
\date{\today}
\dedicatory{To Anna} 

\begin{abstract}
We prove that the Higman-Thompson groups $G_{n,r}^+$ and $G_{m,s}^+$ are isomorphic if and only if $m=n$ and $\mbox{gcd}(n-1,r)=\mbox{gcd}(n-1,s)$.
\end{abstract}

\maketitle

\section*{Introduction}

Finitely presented groups are one of the most important classes of infinite groups, both by its ubiquity (e.g. they are fundamental groups of compact manifolds) and by the number of interesting subclasses it contains (hyperbolic groups and automatic groups, among others). Also they have interest in connection with algorithmic properties, as showed by Higman \cite{HigmanRec}, who stated that a finitely generated group is embeddable in a finitely presented group if and only if is recursively presented. As simple groups are one of the milestones in the development of Group Theory, a fundamental topic in the study of groups is that of finitely presented simple groups.

The study of finitely presented simple groups began with Thompson's discovering in 1965 of the firsts two infinite examples in this class \cite{Thompson}, now known as $G_{2,1}$ and $T_{2,1}$. In 1974 Higman \cite{Higman} constructed a countably infinite family of finitely presented simple groups generalizing Thompson's group $G_{2,1}$. These are the commutator subgroups $G_{n,r}^+$ of the groups $G_{n,r}$ introduced in the same paper. There are various ways of describing these groups: automorphism groups of $r$-generated free algebras in the variety of algebras of sets that are in bijection with its own $n$-th direct power \cite{Higman}, groups of piecewise linear homeomorphisms of the unit interval with prescribed slopes and limited sets of non-differenciable points \cite{Stein}, groups of tree diagrams of finite $n$-ary $r$-forests \cite{CFP}, or groups of maximal inescapable (cofinite) isomorphisms \cite{Scott}.

Unfortunately, these groups are still not fully understood. For example, the isomorphism problem is not completely solved in this class, even if Higman \cite{Higman} highlighted a great part of it. Higman showed that showed that this family contains infinitely many isomorphism types.  Also he showed that $G_{n,r}^+\cong G_{m,s}^+$ implies that $m=n$ and $\mbox{gcd}(r,n-1)=\mbox{gcd}(s,n-1)$ \cite[Theorem 6.4]{Higman}, while the converse is known only for some particular cases (e.g. when $r\equiv s$ (mod $n-1$) \cite[Section 3]{Higman}, or when $s=rc$ with $c$ a divisor of $n$ \cite[Theorem 7.3]{Higman}). 

In this paper we prove that the converse of \cite[Theorem 6.4]{Higman} holds: $G_{n,r}^+$ and $G_{m,s}^+$ are isomorphic if and only if $m=n$ and $\mbox{gcd}(n-1,r)=\mbox{gcd}(n-1,s)$. Hence, we close the isomorphism problem for this class. The key point for proving this result relies in the connection of this problem with a longstanding problem about isomorphisms of finitely presented algebras stated by Leavitt \cite{L1,L2}, and recently solved by Abrams, \'Anh and the author \cite{AbAnhP}.

Now we summarize the contents of this paper. In Section 1 we recall the definition of $G_{n,r}$ following the lines of \cite{Scott}, and we list some properties enjoyed by these groups. In Section 2 we recall the definition of Leavitt algebras, and we quote  \cite[Theorem 4.14]{AbAnhP}:
\begin{quotation}
Let $d,n$ be positive integers, and $K$ any field.  Let $L_{K,n} = L_n$ denote the Leavitt algebra of type $(1,n-1)$
with coefficients in $K$. Then $L_n\cong M_d(L_n)$ if and only if $\mbox{gcd}(d,n-1)=1$.
\end{quotation}
Since the isomorphism is explicitly given in terms of the generators of the algebra, we use it in Section 3, where we relate $G_{n,r}$ with a group of invertible matrices in $M_r(L_n)$, and as a byproduct we give a proof of the converse of \cite[Theorem 6.4]{Higman}.

\section{Basics on Higman-Thompson groups}

We will fix the essential definitions and results about Higman-Thompson groups that we will need in the sequel. Our sources are \cite{CFP, Higman, Rover, Scott}.

Let $n, r \in \N$, $n\geq 2$ and $r\geq 1$, let $\mathcal{A}_n=\{ a_1, \dots ,a_n\}$ be an alphabet (through the rest of the paper we will assume $\mathcal{A}_n=\{ 1, \dots ,n\}$  by defect), let $X_r=\{ x_1, \dots ,x_r\}$ a set of $r$ elements disjoint of $\mathcal{A}_n$, and let $W_n$ be the free monoid generated by $\mathcal{A}_n$ and the empty word. Denote by $X_rW_n$ the set of finite words of the form $x_i\alpha$, where $\alpha\in W_n$. 

Given $u,v\in X_rW_n$, we will denote $u\leq v$ if there exists $\alpha \in W_n$ such that $v=u\alpha$; notice that $\leq$ is a partial order, and so $u<v$ means $u\leq v$ and $u\ne v$ with no ambiguity. We will say that a subset $B$ of $X_rW_n$ is independent if its elements are pairwise $\leq$-incomparable.

A nonempty subset $V$ of $X_rW_n$ is said to be a subspace if it is closed under right multiplication by elements of $W_n$. A subset $B$ of a subspace $V$ is a basis if it is independent and $V=BW_n$; a subset $B$ of $X_rW_n$ is a basis if there exists a subspace $V$ for which $B$ is a basis. Notice that the set $B_V=\{ y\in V\mid \mbox{ no proper initial segment of  }y \mbox{ belongs to } V\}$ is a basis for $V$, so that every subspace has a basis. In particular, $X_r$ is a basis for $X_rW_n$.

A subspace $V$ is cofinite if $\vert X_rW_n\setminus V\vert<\infty$. A basis $B$ is cofinite if $V=BW_n$ is a cofinite subspace. Notice that then $V$ is cofinite if it has a maximal finite basis. In particular, any finite basis is contained in a cofinite basis. Also it is clear that any finite intersection of cofinite subspaces is a cofinite subspace (this is \cite[Corollary 1 to Lemma 2.4]{Higman} stated in a different language).

If $u$ is an element of $X_rW_n$, we will say that $\{ua_1, \dots ,ua_n\}$ is a simple expansion of $u$. Given a basis $B$ and $u\in B$ any element, $B'=(B\setminus \{ u\})\cup \{ua_1, \dots ,ua_n\}$ is again a basis, that we call a simple expansion of $B$. Given $B,C$ basis, we say that $C$ is an expansion of $B$ if there is a finite chain $B_0, \dots ,B_k$ of basis such that $B_0=B$, $B_k=C$ and $B_{i+1}$ is a simple expansion of  $B_i$ for every $0\leq i\leq k-1$. Of course, any expansion of a cofinite basis is a cofinite basis as well.

An homomorphism $\theta$ between subspaces of $X_rW_n$ is a map satisfying $(uw)\theta=(u\theta)w$ for all $w\in W_n$, whenever $u\theta$ is defined. An isomorphism is a bijective homomorphism, and if the domain and the range of an isomorphism is cofinite, we say that it is a cofinite isomorphism. An extension of a cofinite isomorphism $\theta$ is a cofinite isomorphism $\theta'$ such that $u\theta'=u\theta$, whenever $u\theta$ is defined. A cofinite isomorphism is maximal if it has no nontrivial extensions.

Now, we quote two fundamental facts:

\begin{lem}\label{L:Scott1} {\rm (c.f. \cite[Lemma 1]{Scott})} Every cofinal isomorphism $\theta$ has a unique maximal extension $\theta^*$.
\end{lem}

Let $\phi_i: U_i\rightarrow V_i$ be cofinite isomorphisms for $i=\{ 1,2\}$. Now fix $S=V_1\cap U_2$, $R=S\phi_1^{-1}$, $T=S\phi_2$, which are cofinite subspaces, and notice that $({\phi_1}_{\vert R})\circ ({\phi_2}_{\vert S})$ is a cofinal isomorphism from $R$ to $S$. So, we define $\phi_1\phi_2:=(({\phi_1}_{\vert R})\circ ({\phi_2}_{\vert S}))^*$. With this definition we have

\begin{lem}\label{L:Scott2} {\rm (c.f. \cite[Lemma 2]{Scott})} The set of maximal cofinite isomorphisms is a group under the above defined operation.
\end{lem}

The group defined in Lemma \ref{L:Scott2} is the Higman-Thompson group $G_{n,r}$ defined originally in \cite{Higman}. We introduce a the representation of the elements of $G_{n,r}$ which turns out to be a useful instrument to deal with the group. 

Whenever $B=\{ y_1, \dots ,y_N\}$ and  $C=\{ z_1, \dots ,z_N\}$ are expansions of $X_r$ (and thus cofinite basis), the bijection
$$
\begin{array}{cccc}
 \theta &B  & \rightarrow  & C  \\
 & y_i & \mapsto  & z_i   
\end{array}
$$
extends naturally to a cofinite isomorphism $\theta: BW_n\rightarrow CW_n$, so that $\theta \in G_{n,r}$. Thus, we can represent $\theta$ by the symbol
$$\theta =\left(
\begin{array}{ccc}
 y_1 & \dots  & y_N  \\
 z_1 & \dots  & z_N  
\end{array}
\right).
$$
Conversely, every element $\theta\in G_{n,r}$ admits such a representation \cite[Lemma 4.1]{Higman}. 

Whenever 
$$\varphi =\left(
\begin{array}{ccc}
 x_1 & \dots  & x_M  \\
 t_1 & \dots  & t_M  
\end{array}
\right)
$$
is a symbol for any other element in $G_{n,r}$, \cite[Corollary 1 to Lemma 2.4]{Higman} guarantees that there exists a common expansion $\{ s_1 , \dots  ,s_P\} $ of $\{ z_1 , \dots  ,z_N\} $ and $\{ x_1 , \dots  , x_M\} $ so that
$$\theta =\left(
\begin{array}{ccc}
 y'_1 & \dots  & y'_P  \\
 s_1 & \dots  & s_P  
\end{array}
\right)
\mbox{ and }
\varphi =\left(
\begin{array}{ccc}
 s_1 & \dots  & s_P  \\
 t'_1 & \dots  & t'_P  
\end{array}
\right)
$$
and thus
$$\theta \varphi=\left(
\begin{array}{ccc}
 y'_1 & \dots  & y'_P  \\
  t'_1 & \dots  & t'_P  
\end{array}
\right).
$$

A relevant subgroup of $G_{n,r}$ is the commutator subgroup, usually denoted by $G_{n,r}^+$. At it was shown in \cite{Higman} (c.f. \cite[Lemma 2.1]{Rover}), the index of $G_{n,r}^+$ in $G_{n,r}$ is $\mbox{gcd}(n-1, 2)$, so that $G_{n,r}^+$ coincides with $G_{n,r}$ whenever $n$ is even. For the sake of uniform notation (c.f. \cite[Section 5]{Higman}), we write $G_{n,r}^+=G_{n,r}$ when $n$ is even.

For any $n\geq 2, r\geq 1$, some interesting features enjoyed by these groups are the following:
\begin{enumerate}
\item The group $G_{n,r}$ is finitely presented \cite[Theorem 4.6]{Higman}.
\item The group $G_{n,r}^+$ is simple \cite[Theorem 5.4]{Higman}.
\item The group $G_{n,r}^+$ contains an isomorphic copy of every countable locally finite group \cite[Theorem 6.6]{Higman}.
\item The defining relations of the group $G_{n,r}$ are recursively enumerable, so that $G_{n,r}$ has soluble word problem, and thus conjugacy and order soluble problems \cite[Section 9]{Higman}.
\end{enumerate}

\section{Isomorphisms of Leavitt algebras}

We begin by defining the Leavitt algebras $L_{K}(1,n)$, which was investigated originally by Leavitt in his seminal paper \cite{L1}.
For any positive integer $n\geq 2$, and field $K$, we denote
$L_K(1,n)$ by $L_{K,n}$, and call it the Leavitt algebra of
type (1,n-1) with coefficients in $K$. (When $K$ is understood,
we denote this algebra simply by $L_n$).  Precisely, $L_{K,n}$ is
the quotient of the free associative $K$-algebra in $2n$
variables:
$$L_{K,n}=K<X_1,...,X_n,Y_1,...,Y_n>/T,$$
where $T$ is the ideal generated by the relations $ X_iY_j - \delta_{ij}1_K$ (for $1\leq i,j \leq n$) and
$\sum_{j=1}^{n} Y_jX_j - 1_K$. The images of $X_i,Y_i$ in $L_{K,n}$ are denoted respectively by $x_i,y_i$. In
particular, we have the equalities $x_iy_j = \delta_{ij}1_K$ and $\sum_{j=1}^{n} y_jx_j = 1_K$ in $L_n$. A multiindex will be a sequence $I=\{ i_1\dots, i_k\}$ with $i_j\in\{ 1, 2, \dots ,n\}$ for all $1\leq j\leq k$. We will then denote $y_I=y_{i_1}y_{i_2}\cdots y_{i_k}$ and $x_I=x_{i_k}x_{i_{k-1}}\cdots x_{i_1}$.

We now fix a fundamental property of $L_n$ that is basic for our purposes.

\begin{lem}\label{Basicprops}  {\rm (\cite[Theorem 8]{L1})} Let $K$ be any field, let $n\geq 2$ be a natural number. Then, $L_n$ has module type $(1,n-1)$.  In particular, if $r\equiv s$ (mod $n-1$) then $L_n^r \cong
L_n^s$ as free left $L_n$-modules.  Consequently, if $r\equiv s$ (mod $n-1$), then there is an isomorphism of matrix
rings $M_r(L_n)\cong M_s(L_n)$.
\end{lem}

\begin{rema}\label{R:special iso}
{\rm Suppose that $s=r+(n-1)$, and denote $\widehat{x}=(x_1,\dots ,x_n)$ and $\widehat{y}=(y_1, \dots ,y_n)$. Then, abovementioned isomorphism is given by the rule
$$
\begin{array}{cccc}
 \varphi: & M_r(L_n)  & \rightarrow  & M_s(L_n)  \\
 & A & \mapsto  &   \mbox{diag}(I_{r-1}, \widehat{x}^t)\cdot A\cdot \mbox{diag}(I_{r-1}, \widehat{y})
\end{array}.
$$
In particular, whenever the entries of $A$ has the form $\sum\limits_{i=1}^{k}y_{I_i}x_{J_i}$ (for $\{ I_i, J_i\}_{1\leq i\leq k}$ sets of multiindices), then so are the entries of $\varphi (A)$.
By recurrence on this argument, the same consequence holds for any pair of natural numbers $r,s$ such that  $r\equiv s$ (mod $n-1$). }
\end{rema}

\begin{defi}\label{involutiondef}  {\rm For any field $K$, the extension of the assignments $x_i \mapsto y_i=x_i^*$ and $y_i \mapsto x_i=y_i^*$
 for $1\leq i \leq n$ yields an involution $*$ on $L_K(1,n)$. This involution on $L_K(1,n)$ produces an involution
 on any sized matrix ring $M_m(L_K(1,n))$ over $L_K(1,n)$ by setting $X^* =
(x_{j,i}^*)$  for each $X=(x_{i,j})\in M_m(L_K(1,n))$. We note that if $K$ is a field with involution (which we also denote by $*$), then a second involution on $L_K(1,n)$
may be defined by extending the assignments $k\mapsto k^*$ for all $k\in K$,  $x_i \mapsto y_i=x_i^*$ and $y_i \mapsto
x_i=y_i^*$ for $1\leq i\leq n$. We will say that $X\in M_d(L_n)$ is a unitary provided that $XX^*=X^*X=I_d$, and we will denote by $U(M_d(L_n))$ the group of unitaries of $M_d(L_n)$.}
\end{defi}

Now, we will quote the essential result of this section

\begin{theor}\label{Th:IsoLeavitt} {\rm (\cite[Theorem 4.14]{AbAnhP})}
Let $d,n$ be positive integers, and $K$ any field.  Let $L_{K,n} = L_n$ denote the Leavitt algebra of type $(1,n-1)$
with coefficients in $K$. Then $L_n\cong M_d(L_n)$ if and only if $\mbox{gcd}(d,n-1)=1$.
\end{theor}

Let us fix the details needed to prove Theorem \ref{Th:IsoLeavitt}, in order to explain why it is so important in the proof of our main result. Essentially, we need to construct a $K$-algebra isomorphism
$$
\begin{array}{cccc}
 \varphi: &  L_n & \rightarrow & M_d(L_n)   \\
 & x_i & \mapsto  &   X_i\\
 & y_j & \mapsto  &   Y_j
\end{array}.
$$
Since $L_n$ is a simple algebra, it is enough to fix a set $\{ X_1, \dots, X_n,Y_1,\dots ,Y_n\}\subset M_d(L_n)$ satisfying the definitory relations of the generators of $L_n$, and generating $M_d(L_n)$. Now, we present the appropriate $2n$ matrices. For any unital ring $R$ and $i\in\{1,2,...,d\}$ we denote the
idempotent $e_{i,i}$  of the matrix ring $M_d(R)$ simply by
$e_i$. We write
$n=qd+r$ with $2\leq r \leq d$.  We assume $d<n$, so that $q\geq 1$. The matrices $X_1, X_2, ..., X_q$ are given as
follows. For $1\leq i \leq q$ we define
$$X_i=
\begin{pmatrix}x_{(i-1)d+1}&0& &0 \\
        x_{(i-1)d+2}&0& &0 \\
        \vdots&0&...&0 \\
        x_{id}&0& &0
        \end{pmatrix}
= \sum_{j=1}^{d}x_{(i-1)d+j}e_{j,1} $$ The two matrices $X_{q+1}$
and $X_{q+2}$ play a pivotal role here.  They are defined as
follows.
$$X_{q+1}=
\begin{pmatrix}x_{qd+1}&0&0& &0&0& &0& \\
        x_{qd+2}&0&0& &0&0& &0& \\
        \vdots&0&0& &0&0& &0& \\
        x_n&0&0&...&0&0&...&0&\\
        0&1&0& &0&0& &0& \\
        0&0&1& &0&0& &0& \\
        & & &\vdots& & & & & \\
        0&0&0& ...&1&0& &0& \end{pmatrix}$$
$$ = \sum_{i=1}^{d-r}e_{i+r,i+1} + \sum_{t=1}^{r}x_{qd+t}e_{t,1}$$
and
$$X_{q+2}=
\begin{pmatrix}0& &0&1&0&0& &0&0 \\
        0& &0&0&1&0& &0&0 \\
         & & & & &\vdots& & &\\
        0& &0&0&0&0& &1&0 \\
        0&...&0&0&0&0& &0&a_{q+2,r-1} \\
        0& &0&0&0&0& &0&a_{q+2,r} \\
         & & &\vdots& & & & &\vdots \\
        0& &0&0&0&0& &0&a_{q+2,d}\end{pmatrix}$$
$$ = \sum_{j=1}^{r-2}e_{j,j+s} + \sum_{t=1}^{d-(r-2)}a_{q+2,(r-2)+t}e_{(r-2)+t,d}$$
(where the elements $a_{q+2,r-1}, a_{q+2,r},..., a_{q+2,d} \in
L_n$ are monomials in $x$-variables).  In case $d-r=0$ or $r-2=0$ we interpret the appropriate
sums as zero.

The remaining matrices $X_{q+3},...,X_n$ will have the same general
form. In particular, for $q+3 \leq i \leq n$,
$$X_i=
\begin{pmatrix}0& &0&a_{i,1} \\
        0& &0&a_{i,2} \\
         0&...&\vdots&\\
        0& &0&a_{i,d}
        \end{pmatrix}
= \sum_{j=1}^{d}a_{i,j}e_{j,d} $$ (where the elements $a_{i,1},
a_{i,2}, ... ,a_{i,d} \in L_n$ are monomials in the $x$-variables).  In case $q+3 > n$ then we
understand that there are no matrices of this latter form in our
set of $2n$ matrices.  We note that we always have the matrices
$X_{q+1}$ and $X_{q+2}$, since $n=qd+r\geq q\cdot 1 + 2$. We define the matrices $Y_i$ for $1\leq i \leq n$ by setting $Y_i
= X_i^*$. 

Now consider this set, which we will call ``The List":

$$x_1^{d-1}$$
$$x_2x_1^{d-2}, x_3x_1^{d-2}, ... , x_nx_1^{d-2}$$
$$x_2x_1^{d-3}, x_3x_1^{d-3}, ... , x_nx_1^{d-3}$$
$$\vdots$$
$$x_2x_1, x_3x_1, ... , x_nx_1$$
$$x_2, x_3, ... , x_n$$

The key of the proof of Theorem \ref{Th:IsoLeavitt} is that, whenever $\mbox{gcd}(d,n-1)=1$, there is a rule to assign an element in The List to each $a_{i,j}$ in the above set of matrices, in such a way that the resulting set $\{ X_1, \dots, X_n,Y_1,\dots ,Y_n\}$ satisfies the definitory relations of the generators of $L_n$, and generates $M_d(L_n)$. Thus, under such a choice, the above defined map $\varphi$ is a $K$-algebra isomorphism. 

\begin{rema}\label{R:special iso2}
{\rm Because of the definition of the above mentioned isomorphism, it is clear that whenever $a\in L_n$ has the form $\sum\limits_{i=1}^{k}y_{I_i}x_{J_i}$ (for $\{ I_i, J_i\}_{1\leq i\leq k}$ sets of multiindices), then the entries of $\varphi (a)$ have the same form. }
\end{rema}

We will prove an easy consequence of Theorem \ref{Th:IsoLeavitt} that will be useful in the sequel. For, we quote the following fact

\begin{lem}\label{L:AbGroups} {\rm (\cite[Lemma 1]{AbSmith})} Let $G$ be a finitely generated abelian group (written additively). Let $x\in G$ be an element of finite order $n$, and let $c,d\in \N$. There exists an automorphism $\varphi: G\rightarrow G$ with $\varphi (cx)=dx$ if and only if $\mbox{gcd}(c,n)=\mbox{gcd}(d,n)$.
\end{lem}

\begin{corol}\label{Corol:IsoLeavitt}
Let $n, r, s$ be positive integers, and $K$ any field.  Let $L_{K,n} = L_n$ denote the Leavitt algebra of type $(1,n-1)$
with coefficients in $K$. If $\mbox{gcd}(r,n-1)=\mbox{gcd}(s,n-1)$, then $M_r(L_n)\cong M_s(L_n)$.
\end{corol}
\begin{proof}
By Lemma \ref{L:AbGroups}, applied to $G=\Z/(n-1)\Z$, $x=[1]\in \Z/(n-1)\Z$, $c=r$ and $d=s$, there exists a group automorphism $\varphi: \Z/(n-1)\Z\rightarrow \Z/(n-1)\Z$ such that $\varphi ([r])=[s]$. Thus, there exists $l\in \N$ with $\mbox{gcd}(l,n-1)=1$ such that $[lr]=[s]$. Since $lr\equiv s$ (mod $n-1$), we have $M_s(L_n)\cong M_r(M_l(L_n))$ by Lemma \ref{Basicprops}. Now, $L_n\cong M_l(L_n)$ by Theorem \ref{Th:IsoLeavitt}, so that $M_r(M_l(L_n))\cong M_r(L_n)$, which completes the proof.
\end{proof}

\begin{rema}\label{R:special iso3}
{\rm Because of Remarks \ref{R:special iso} and \ref{R:special iso2}, the isomorphism given by Corollary \ref{Corol:IsoLeavitt} has the property that, whenever the entries of $A$ has the form $\sum\limits_{i=1}^{k}y_{I_i}x_{J_i}$ (for $\{ I_i, J_i\}_{1\leq i\leq k}$ sets of multiindices), then so are the entries of $\varphi (A)$. This fact play a role in the proof of the main result of this paper.}
\end{rema}

\section{The main result}

In this section, we will prove the main result of the paper.

\begin{defi}\label{D:grupunitaris}
{\rm Let $n\geq 2, r\geq 1$ be natural numbers. We denote by $\mathcal{P}_{n,r}$ the subset of the group $U(M_r(L_n))$ of unitaries of $M_r(L_n)$ composed by matrices in which all the entries are either $0$ or have the form 
$$\sum\limits_{i=1}^{m}y_{I_i}x_{J_i},$$
where the $I_i, J_i$ are multiindices.}
\end{defi}

\begin{lem}\label{L:grupunitaris}
For any $n\geq 2, r\geq 1$ natural numbers, $\mathcal{P}_{n,r}$ is a subgroup of $U(M_r(L_n))$.
\end{lem}
\begin{proof}
Since $1=\sum\limits_{i=1}^{n}y_ix_i$, it is clear that the identity matrix $I_r$ belongs to $\mathcal{P}_{n,r}$, whence it is a nonempty set.

Fix $X,Y\in \mathcal{P}_{n,r}$ two elements. Since $Y$ is an unitary, $Y^{-1}$ is the conjugated transpose of $Y$ (so it lies in $\mathcal{P}_{n,r}$ too), and hence the entries in $XY^{-1}$ are of the form $\sum\limits_{k=1}^{r}a_{i,k}b_{k,j}$, where
$$a_{i,k}b_{k,j}=\left( \sum\limits_{i=1}^{m}y_{I_i}x_{J_i}\right)\cdot \left( \sum\limits_{i=1}^{m'}y_{I'_i}x_{J'_i}\right).$$
As
$$(\ast )\hspace{.2truecm} x_{J_r}y_{I'_s}=\left\{\begin{array}{lll}  y_{\widehat{I}_s} & \mbox{if} &
I'_s=J_r\widehat{I}_s\\
 x_{\widehat{J}_r} & \mbox{if} &
J_r=\widehat{J}_rI'_s\\
 0 &  & \mbox{otherwise}\\
\end{array}
\right.
$$
we conclude that $XY^{-1}\in \mathcal{P}_{n,r}$, as desired.
\end{proof}

Notice that, if we fix the alphabet $\mathcal{A}_n=\{ 1, \dots ,n\}$, then each multiindex is an element of $W_n$, and the identity $(\ast )$ in Lemma \ref{L:grupunitaris} says that $x_{J_r}y_{I'_s}=0$ if and only if $x_1I'_s$ and $x_1J_r$ are independent elements of $X_1W_n$. The key for connecting the isomorphism problem of Higman-Thompson groups with Leavitt algebras lies precisely in this fact, that we will exploit.

Now, we will prove a technical results that will be needed later.

\begin{lem}\label{L:InverseMap}
Let $n\geq 2$ be a natural number. If $\alpha =\sum\limits_{i=1}^{m}y_{I_i}x_{J_i}\in \mathcal{P}_{n,1}$, then both $\{ I_1, \dots , I_m\}$ and $\{ J_1, \dots , J_m\}$ are expansions of the basis $\{ x_1\}$ of $X_1W_{n}$ (and thus basis).
\end{lem}
\begin{proof}
As the argument is symmetric, we will proof it only for $\{ I_1, \dots , I_m\}$.

First, suppose that $\{ I_1, \dots , I_m\}$ do not contain a complete expansion of $\{ x_1\}$. Then, two different cases could happen:
\begin{enumerate}
\item The set $\{ I_1, \dots , I_m\}$ is independent, and thus can be completed to a basis $$\{ I_1, \dots , I_l, \widehat{I}_{l+1}, \dots ,\widehat{I}_{r}\}.$$
\item The set $\{ I_1, \dots , I_m\}$ is not independent. So, we can chose a maximal independent subset $\{ I_1, \dots , I_l\}\subsetneq\{ I_1, \dots , I_m\}$. 
\end{enumerate}
Hence, in any of both cases, for $l\leq m$ there exist a maximal independent subset $\{ I_1, \dots , I_l\}\subseteq\{ I_1, \dots , I_m\}$ and a multiindex $Z$ such that $\{ I_1, \dots , I_l, Z\}$ can be expanded to a basis. But then, as $Z$ is independent of the $I_j$s,
$$0\ne x_Z=x_Z\cdot 1=x_Z(\alpha \alpha^*)=(x_Z\alpha )\alpha^*=0\alpha ^*=0$$
which is impossible.

Now suppose that $\{ I_1, \dots , I_m\}$ contains an expansion (i.e. a basis) but it is not a basis, i.e it is not an independent set. Fix $\{ I_1, \dots , I_l\}\subsetneq\{ I_1, \dots , I_m\}$ a basis, and notice that
$$
1=\left( \sum\limits_{i=1}^{l}y_{I_i}x_{J_i}\right)\cdot \left( \sum\limits_{j=1}^{l}y_{J_j}x_{I_j}\right).
$$
Hence,
\begin{align*}
1= &\alpha\alpha^*=  \left( \sum\limits_{i=1}^{m}y_{I_i}x_{J_i}\right)\cdot \left( \sum\limits_{j=1}^{m}y_{J_j}x_{I_j}\right)\\& =\left( \sum\limits_{i=1}^{l}y_{I_i}x_{J_i}+\sum\limits_{i=l+1}^{m}y_{I_i}x_{J_i}\right)\cdot \left( \sum\limits_{j=1}^{l}y_{J_j}x_{I_j}+\sum\limits_{j=l+1}^{m}y_{J_j}x_{I_j}\right)\\& = 1+ \sum\limits_{i=1}^{l}\sum\limits_{j=l+1}^{m}y_{I_i}x_{J_i}y_{J_j}x_{I_j} + \sum\limits_{i=l+1}^{m}\sum\limits_{j=1}^{l}y_{I_i}x_{J_i}y_{J_j}x_{I_j} + \sum\limits_{i=l+1}^{m}\sum\limits_{j=l+1}^{m}y_{I_i}x_{J_i}y_{J_j}x_{I_j}.
\end{align*}
Thus, the last 3 summands equal zero, and in particular for any $l+1\leq i\leq m$ we have $0=y_{I_i}x_{J_i}y_{J_i}x_{I_i}=y_{I_i}x_{I_i}$, which is impossible. So, we are done.
\end{proof}

The goal is to prove that for any $n\geq 2, r\geq 1$, $\mathcal{P}_{n,r}\cong G_{n,r}$. In order to do more comprehensible the argument, first we will prove the result in the particular case $r=1$. This result is analogous to \cite[Proposition 9.6]{N}, but the proof is different.

\begin{prop}\label{P:IsoFor r=1}
If $n\geq 2$ is a natural number, then $\mathcal{P}_{n,1}\cong G_{n,1}$.
\end{prop}
\begin{proof}
By \cite[Lemma 4.1]{Higman}, given an element $x\in G_{n,1}$, we can express it by using a symbol
$$x=
\left(
\begin{array}{ccc}
I_1  &  \dots  & I_m  \\
 J_1 & \dots   & J_m  
\end{array}
\right)
$$
where both $\{ I_1, \dots , I_m\}$ and $\{ J_1, \dots , J_m\}$ are expansions of the basis $\{ x_1\}$ of $X_1W_{n}$ (and thus basis). Now define
$$\alpha_x=\sum\limits_{i=1}^{m}y_{I_i}x_{J_i}\in L_n.$$
Notice that 
$$
\alpha_x\alpha_x^*=\left( \sum\limits_{i=1}^{m}y_{I_i}x_{J_i}\right)\cdot \left( \sum\limits_{j=1}^{m}y_{J_j}x_{I_j}\right)= \sum\limits_{i=1}^{m}y_{I_i}x_{I_i}=1,
$$
where the last two equalities are due to the fact that both $\{ I_1, \dots , I_m\}$ and $\{ J_1, \dots , J_m\}$ are expansions of the basis $\{ x_1\}$ of $W_{n,1}$. Similarly, $\alpha_x^*\alpha_x=1$, so that $\alpha_x\in \mathcal{P}_{n,1}$. Define a map
$$\begin{array}{cccc}
 \varphi: &  G_{n,1} & \rightarrow & \mathcal{P}_{n,1}   \\
  &  x & \mapsto  & \alpha_x  
\end{array},
$$
and notice that $\varphi$ send symbols equivalent by elementary expansions to the same element in $L_n$. Thus, if 
$$x=
\left(
\begin{array}{ccc}
I_1  &  \dots  & I_m  \\
 J_1 & \dots   & J_m  
\end{array}
\right)
\mbox{ and }
y=
\left(
\begin{array}{ccc}
R_1  &  \dots  & R_k  \\
 S_1 & \dots   & S_k  
\end{array}
\right),
$$
again by \cite[Corollary 1 to Lemma 2.4]{Higman} there exists a common expansion $\{ J'_1 , \dots  , J'_t\}$ of both $\{  J_1, \dots   , J_m  \}$ and $\{ R_1  ,  \dots  , R_k\}$ such that
$$x=
\left(
\begin{array}{ccc}
I'_1  &  \dots  & I'_t  \\
 J'_1 & \dots   & J'_t  
\end{array}
\right)
\mbox{ and }
y=
\left(
\begin{array}{ccc}
J'_1 & \dots   & J'_t  \\
 S'_1 & \dots   & S'_t  
\end{array}
\right),
$$
whence 
$$xy=
\left(
\begin{array}{ccc}
I'_1  &  \dots  & I'_t  \\
  S'_1 & \dots   & S'_t  
\end{array}
\right).
$$
By the above remark, we get $\varphi (xy)=\varphi (x)\varphi (y)$, so that $\varphi$ is a group morphism.

Now, if $\alpha =\sum\limits_{i=1}^{m}y_{I_i}x_{J_i}\in \mathcal{P}_{n,1}$, then the element
$$x_{\alpha}=
\left(
\begin{array}{ccc}
I_1  &  \dots  & I_m  \\
 J_1 & \dots   & J_m  
\end{array}
\right)
$$
belong to $G_{n,1}$ by Lemma \ref{L:InverseMap}, so that 
$$\begin{array}{cccc}
 \psi: &  \mathcal{P}_{n,1} & \rightarrow & G_{n,1}   \\
  &  \alpha & \mapsto  & x_{\alpha  }
\end{array}
$$
is a well-defined map. Moreover, $\varphi(x_{\alpha})=\alpha$, so that $\varphi $ is an onto map. As $\psi$ is clearly compatible with the equivalence of symbols by elementary expansions, in turns out that $\psi$ is a group morphism. A simple inspection shows that $\varphi$ and $\psi$ are mutually inverses, so we are done.
\end{proof}

Now, we prove the general version.

\begin{prop}\label{Prop:IsoGroups}
If $n\geq 2$ and $r\geq 1$ are natural numbers, then $\mathcal{P}_{n,r}\cong G_{n,r}$.
\end{prop}
\begin{proof}
Consider $X_r=\{ x_1, \dots ,x_r\}$ as basis of $X_rW_{n}$, and take
$$x=
\left(
\begin{array}{ccc}
I_1  &  \dots  & I_k  \\
 J_1 & \dots   & J_k  
\end{array}
\right)\in G_{n,r}.
$$
For any $1\leq i,j\leq r$, consider the strictly ascending finite sequence
$$1\leq l(i,j)_1<\cdots < l(i,j)_{s(i,j)}\leq k$$
such that $I_{l(i,j)_t}$ starts in $x_i$ and $J_{l(i,j)_t}$ starts in $x_j$ for every $1\leq t\leq s(i,j)$ (Notice that it can happens for some sequences in this list to be empty). Consider now the multiindices $I'_{l(i,j)_t}$ and $J'_{l(i,j)_t}$ obtained from $I_{l(i,j)_t}$ and $J_{l(i,j)_t}$ (respectively) by erasing the initial $x_i$ (resp. $x_j$). Now we define the matrix $X\in M_r(L_n)$ whose $(i,j)$-entry is
$$X_{i,j}=\sum\limits_{p=1}^{s(i,j)}y_{I'_{l(i,j)_p}}x_{J'_{l(i,j)_p}}.$$
Notice that $x_{J'_{l(i,k)_p}}y_{J'_{l(j,k)_q}}=\delta_{i,j}\cdot\delta_{p,q}$; indeed, if $x_{J'_{l(i,k)_p}}y_{J'_{l(j,k)_q}}=1$ for $i\ne j$ or $p\ne q$, then the symbol of $x$ shall contain entries
$$
x=
\left(
\begin{array}{ccccc}
\dots & x_iI'_{l(i,k)_p}  &  \dots  & x_jI'_{l(j,k)_q} &\dots \\
 \dots & x_kJ'_{l(i,k)_p} & \dots   & x_kJ'_{l(j,k)_q}  & \dots
\end{array}
\right)
$$
with $x_kJ'_{l(i,k)_p}=x_kJ'_{l(j,k)_q}$, which is impossible by definition of symbol. Also, since $\{ I_1  ,  \dots  , I_k \}$ is an expansion of $X_r$, 
$$\{ I_{l(i,1)_1}, \dots, I_{l(i,1)_{s(i,1)}},\dots ,I_{l(i,r)_1}, \dots, I_{l(i,r)_{s(i,r)}}\}$$
is an expansion of the element $x_i\in X_r$.

We then have
\begin{align*}
(XX^*)_{i,j}= &\sum\limits_{k=1}^{r}\left( \sum\limits_{p=1}^{s(i,k)} y_{I'_{l(i,k)_p}}x_{J'_{l(i,k)_p}}\right)\cdot \left( \sum\limits_{q=1}^{s(j,k)} y_{J'_{l(j,k)_q}}x_{I'_{l(j,k)_q}}\right)\\ & = \delta_{i,j} \cdot\sum\limits_{k=1}^{r} \sum\limits_{p=1}^{s(i,k)} y_{I'_{l(i,k)_p}}x_{I'_{l(i,k)_p}}=\delta _{i,j},
\end{align*}
and similarly $(X^*X)_{i,j}=\delta _{i,j}$. Hence, $X\in \mathcal{P}_{n,r}$. Thus, 
$$\begin{array}{cccc}
 \varphi: &  G_{n,r} & \rightarrow & \mathcal{P}_{n,r}   \\
  &  x & \mapsto  & X  
\end{array}
$$
is a well-defined map. Clearly $\varphi$ respects the equivalence of symbols by elementary expansions in $G_{n,r}$, so that it is straightforward but tedious to prove that in fact it is a group morphism.

Now, take $X\in \mathcal{P}_{n,r}$. For any $1\leq i,j\leq r$,
$$X_{i,j}=\sum\limits_{p=1}^{s(i,j)}y_{I_{l(i,j)_p}}x_{J_{l(i,j)_p}}$$
for suitable sets of multiindices. We will show that both $W_I=\left\{ \left\{ x_i{I_{l(i,j)_p}}\right\}_{1\leq p\leq s(i,j)} \right\}_{1\leq i,j\leq r}$ and $W_J=\left\{ \left\{ x_j{J_{l(i,j)_p}}\right\}_{1\leq p\leq s(i,j)} \right\}_{1\leq i,j\leq r}$ are expansions of the basis $X_r$. Notice that for any $1\leq i\leq r$,
$$(\ast) \hspace{.2truecm} 1=\sum\limits_{k=1}^{r}X_{i,k}X_{i,k}^*=\sum\limits_{k=1}^{r}\left( \sum\limits_{p,q=1}^{s(i,k)}y_{I_{l(i,k)_p}}x_{J_{l(i,k)_p}}y_{J_{l(i,k)_q}}x_{I_{l(i,k)_q}}\right).$$
Fix any $1\leq i\leq r$, and let $W_I(i)=\left\{ \left\{ x_i{I_{l(i,j)_p}}\right\}_{1\leq p\leq s(i,j)} \right\}_{1\leq j\leq r}$. If it do not contain a complete expansion of $x_i$, then the same argument as in Lemma \ref{L:InverseMap} shows that there exist a maximal independent subset $W'$ of $W_I(i)$ and a multiindex $Z$ such that $W'\cup \{Z\}$ is a part of a basis for $x_i$. Hence, 
$$0\ne x_Z\cdot e_{i,i}=(x_Z\cdot e_{i,i})(XX^*)=((x_Z\cdot e_{i,i}X)X^*$$
has $(k,j)$-entry equal to zero for any $k\ne i$ and any $j$, while
$$(x_Z\cdot e_{i,i}X)_{i,k}=\sum\limits_{p=1}^{s(i,k)}x_Zy_{I_{l(i,k)_p}}x_{J_{l(i,k)_p}}=0,$$
which is impossible. On the other side, if $W_I(i)$ contains an expansion of $\{ x_i\}$ but it is not a basis, again the argument in Lemma \ref{L:InverseMap} and the identity $(\ast)$ give us a contradiction. Thus, $W_I(i)$ is an expansion of $x_i$, and thus so is $W_I$ of $X_r$. Similarly we get that $W_J$ is an expansion of $X_r$. Since both sets has the same cardinality,
$$x_X=
\left(
\begin{array}{ccccc}
 x_1I_{l(1,1)_1} &\dots &x_iI_{l(i,j)_p}& \dots& x_rI_{l(r,r)_{s(r,r)}}  \\
  x_1J_{l(1,1)_1} &\dots  &x_jJ_{l(i,j)_p}& \dots& x_rJ_{l(r,r)_{s(r,r)}}   
\end{array}
\right)
$$
is a symbol of an element of $G_{n,r}$, so that
$$\begin{array}{cccc}
 \psi: &  \mathcal{P}_{n,r} & \rightarrow & G_{n,r}   \\
  &  X & \mapsto  & x_X  
\end{array}
$$
is a well-defined map. Moreover, $\varphi (x_X)=X$, so that $\varphi$ is an onto map. As $\psi$ clearly respects the equivalence of symbols by elementary expansions, $\psi$ is a group morphism, and $\varphi$ and $\psi$ are mutually inverses, so we are done.
\end{proof}

Now, we are ready to prove the main result in the paper.

\begin{theor}\label{Th:Main}
Let $n, m\geq 2$ and $r, s\geq 1$ be natural numbers. Then, $G_{m,r}^+\cong G_{n,s}^+$ if and only if $m=n$ and $\mbox{gcd}(n-1,r)=\mbox{gcd}(n-1,s)$.
\end{theor}
\begin{proof}
The ``only if'' part is \cite[Theorem 6.4]{Higman}. 

Now, assume that $\mbox{gcd}(n-1,r)=\mbox{gcd}(n-1,s)$. Then, by Corollary \ref{Corol:IsoLeavitt}, there exists a $K$-algebra isomorphism $\varphi: M_r(L_n)\rightarrow M_s(L_n)$ that, by Remark \ref{R:special iso3}, restricts to a group isomorphism $\phi: \mathcal{P}_{n,r}\rightarrow \mathcal{P}_{n,s}$. As $G_{n,r}\cong \mathcal{P}_{n,r}$ and $G_{n,s}\cong \mathcal{P}_{n,s}$ by Proposition \ref{Prop:IsoGroups}, we conclude that $G_{n,r}\cong G_{n,s}$, and thus $G_{n,r}^+\cong G_{n,s}^+$, as desired.
\end{proof}


\end{document}